\documentclass[12pt]{article}
\usepackage{amsmath,amsfonts,amssymb,amsthm, dsfont, graphicx}
\usepackage{listings,paralist}
\newcommand{\RM}[1]{\MakeUppercase{\romannumeral #1{.}}}
\newtheorem{proposition}{Proposition}[section]
\newtheorem{theorem}[proposition]{Theorem}
\newtheorem{corollary}[proposition]{Corollary}
\newtheorem{lemma}[proposition]{Lemma}
\newtheorem{remark}[proposition]{Remark}

\newtheorem{example}[proposition]{Example}

\newcommand{\nc}{\newcommand}
\nc{\I}{{\bf 1}}
\nc{\bA}{{\mathbf A}}
\nc{\bG}{{ G}}
\nc{\bT}{{ T}}
\nc{\bL}{{\mathbf L}}
\nc{\bS}{{ S}}
\nc{\bN}{{\mathbf N}}
\nc{\bM}{{\mathbf M}}
\nc{\cB}{{\mathcal B}}
\nc{\cH}{{\mathcal H}}
\nc{\cA}{{\mathcal A}}
\nc{\cG}{{\mathcal G}}
\nc{\cS}{{\mathcal S}}
\nc{\cT}{{\mathcal T}}
\nc{\cM}{{\mathcal M}}
\nc{\cK}{{\mathcal K}}
\nc{\cO}{{ O}}
\nc{\R}{{\mathbb R}}
\nc{\Q}{{\mathbb Q}}
\nc{\N}{{\mathbb N}}
\nc{\Z}{{\mathbb Z}}

\DeclareMathOperator{\pr}{pr}

\nc{\BP}{\mathbb{P}}
\nc{\BE}{\mathbb{E}}
\nc{\BQ}{\mathbb{Q}}

\begin{document}

\title{Palm pairs and the general mass-transport principle}

\author{{\sc Daniel Gentner and G\"unter Last},
\footnote{
Postal address: Institut f\"ur Stochastik, Universit\"at Karlsruhe (TH),
76128 Karlsruhe, Germany. 
Email addresses: daniel.gentner@kit.edu, last@math.uni-karlsruhe.de} }

\maketitle

\begin{abstract} 
We consider a lcsc group $\bG$ acting
on a Borel space $\bS$ and on
an underlying $\sigma$-finite measure space.
Our first main result is a transport formula connecting
the Palm pairs of jointly stationary random measures on $\bS$.
A key (and new) technical result is a measurable disintegration
of the Haar measure on $\bG$ along the orbits.
The second main result is an intrinsic characterization of
the Palm pairs of a $\bG$-invariant random measure. 
We then proceed with deriving a general version of
the  mass-transport principle for possibly non-transitive and non-unimodular
group operations first in a deterministic
and then in its full probabilistic form.

\vspace{0.3cm}
\noindent
{\em Keywords: }{random measure, Palm measure, stationarity, invariance 
, locally compact group, operation,  Haar measure, orbit 
, mass-transport principle}

\noindent
{\em Subclass: }{60D05, 60G55, 60G60}

\end{abstract}

\section{Introduction}

Let $\bG$ be a locally compact second countable Hausdorff 
(in short {\em lcsc}) topological group 
operating on a Borel space $(\bS,\cS)$ (under a rather weak technical assumption).
Consider a $\sigma$-finite measure $M$ on $\bS\times\bS$ 
that is invariant under joint shifts of both arguments.
It is helpful to think of $M(C\times D)$ as an amount
of mass transported from $C\in\cS$ to $D\in\cS$. 
Assume first that the group is unimodular and that the
group action is transitive, i.e.\
that there is only one orbit. If $B\in\mathcal{S}$ has
positive and finite invariant measure, then
\begin{align}\label{mtp}
M(B\times S) = M(S\times B).
\end{align}
This {\em mass-transport principle} \cite{B:L:P:S:99,Be:Sch:01}
plays an important role in the study of percolation on graphs.
H\"aggstr\"om \cite{Ha:97} was the first who has used
it for this purpose. Last and Thorisson \cite{LaTh08}
noticed that \eqref{mtp} can also be seen as a special case
of Neveu's classical {\em exchange formula}  \cite{Neveu}.
The exchange formula is a versatile tool in the theory 
of random measures and point processes \cite{Ba:Br:03}.
A general {\em lcsc} group admits a modular function
$\Delta:\bG\rightarrow(0,\infty)$ satisfying \eqref{modular}
below. Still assuming the group action to be 
transitive, \eqref{mtp} generalizes to
\begin{align}\label{mtp2}
\int \tilde\Delta(s,t)\I_B(s)M(d(s,t)) = \int \I_B(t)M(d(s,t)).
\end{align}
where $\tilde\Delta(s,t):=\Delta(g)$ if $gs=t$. This formula can 
be derived from Theorem 4.4 in \cite{La08}.
One purpose of this paper is to establish 
the mass-transport principle \eqref{mtp2} in the general, 
possibly non-transitive and non-unimodular case. 
This principle holds for all $B$ satisfying a natural
symmetry condition and with a suitably generalized
definition of the function $\tilde\Delta(s,t)$, see Theorem \ref{shortMTP}.
Equation \eqref{mtp2} is just a special case
of a more general transport formula for stationary random measures on $S$,
see Theorems \ref{thlastcycle} and \ref{MTP}.
The main aim of this paper is to derive these and
a series of related results in our general setting described
above.

In probability theory stationarity
refers to invariance of the distribution under the shifts induced by $G$.
In this paper we will express stationarity by assuming that
$\bG$ operates measurably on an underlying  
$\sigma$-finite measure space $(\Omega,\mathcal{A},\BP)$,
where $\BP$ is invariant under $\bG$. (We will use a probabilistic
language even though $\BP$ is not assumed to have the finite
total mass $1$.) A {\em stationary random measure} $\xi$ is then just
a $\sigma$-finite kernel from $\Omega$ to $\bS$ which is 
{\em invariant} under joint shifts. We will use this
terminology also for other kernels.
Fundamental objects associated with an invariant random measure $\xi$
are its Palm pairs $(\nu,Q)$, where $\nu$ is a {\em supporting measure}
of $\xi$ (a $\sigma$-finite measure on $S$ equivalent
to $\BE_\BP\xi(\cdot)$) and $Q$ is an appropriate kernel from $S$
to $\Omega$ disintegrating the Campbell measure of
$\xi$, see \eqref{Palm}. Given $\nu$,
such a kernel exists under weak technical assumptions,
see Kallenberg \cite{Kall07}.
In the transitive situation ($Gs=S$ for all $s\in S$) 
the Palm kernel can be obtained from a single measure
on $(\Omega,\mathcal{A})$, the {\em Palm measure} of $\xi$,
by suitable shifts, see e.g.\ \cite{tortrat} (treating the
case of a group acting on itself) and \cite{RoZae90,Kall07,La08b}. 
The seminal paper Mecke \cite{Mecke} is dealing
with an Abelian group operating on itself. This case is
of particular relevance for applications of Palm theory,
see e.g.\  \cite{DVJ,SW08}.
A very general approach to Palm measures and their
invariance properties is taken in Kallenberg \cite{Kall07}.
Much of the notation and terminology of the present
paper stems from this source.

The paper starts with a brief repetition of the basic terminology
for invariant measures and kernels (Section \ref{secterm}) and
for random measures and their Palm pairs (Section \ref{secPalm}).
A key technical result is a measurable and invariant disintegration
of the Haar measure on $\bG$ along the orbits in $\bS$, Theorem \ref{lambdakern}.
The first main result is a {\em transport formula} connecting
the Palm pairs of two jointly stationary random measures $\xi$ 
and $\eta$ on $\bS$,
Theorem \ref{thlastcycle}. This result extends
Theorem 3.6 in \cite{LaTh08} to the more
general case studied in this paper. Corollary \ref{NeveuFormula}
generalizes  Neveu's well-known  
{\em exchange formula}  (see e.g.\ \cite{Neveu}). 
Our second main result is an intrinsic characterization of
Palm pairs of an invariant random measure, Theorem \ref{charMecke}.
This extends Mecke's \cite{Mecke} famous
characterization of Palm measures of stationary random
measures on an Abelian group (cf.\ also \cite{La08} for the
case of a general group and \cite{RoZae90}, \cite{La08b} for
the case of a homogeneous space). 
We then proceed in Section \ref{sectransport}
with deriving a general version of the  {\em mass-transport principle} 
\eqref{mtp2} first in a deterministic form (Subsection \ref{subdet})
and then in its full probabilistic form (Subsection \ref{subtrpri}).
Theorem \ref{MTP} is the third main result of this paper.
Applications of this principle in non-transitive settings 
(partially stationary and isotropic tesselations)
will be provided in the paper \cite{GeLa09}, which is in preparation.

Finally in this introduction we state our basic notation
for measures and kernels.
Let $(T,\cT)$ denote a measurable space. If $\mu$ is a measure
on $(T,\cT)$ and $f:T\rightarrow[-\infty,\infty]$ is measurable
then we denote the integral $\int fd\mu$ by $\mu f\equiv \mu(f)$
whenever it is well-defined. We denote
by $\cT_+$ the space of $\cT$-measurable
$[0,\infty]$-valued functions. For $f\in\cT_+$ we write $f\cdot\mu$
for the measure $A\mapsto\int\I_A(t)f(t)\mu(dt)$. If $\nu$ is another measure
on $T$, then $\mu\sim\nu$ means that both $\mu\ll\nu$ and $\nu\ll\mu$.
If $(T',\cT')$ is another measurable space then a 
{\em kernel} $\kappa$ from $T$ to $T'$ is a 
map $\kappa:\bT\times\cT'\rightarrow[0,\infty]$ such that 
for each $B\in\cT'$ the map $t\mapsto\kappa(t,B)$ is measurable 
and that for each $t\in T$ the setfunction $\mu(t,\cdot)$ 
is a measure on $T'$. A kernel 
from $\bT$ to $\bT$ is usually refered to as a kernel {\em on} $\bT$.
The kernel $K$ is {\em $\sigma$-finite} if for 
each $t\in\bT$  the measure $K(t,\cdot)$ is $\sigma$-finite. In this paper
all kernels are assumed to be $\sigma$-finite.
Let $K$ be a kernel from $(T,\cT)$ to another measurable space
$(T',\cT')$. If $\mu$ is a measure on $(T,\cT)$
then $\mu\otimes K$ denotes the measure on the product space 
$(T\times T',\cT\otimes \cT')$ defined by
$\mu\otimes K(A)=\iint \I_A(s,t)K(s,dt)\mu(ds)$, where
$\I_A$ is the indicator function of $A\in \cT\otimes \cT'$.
(Note that $\sigma$-finiteness of $K$ implies measurability of 
$s\mapsto\int \I_A(s,t)K(s,dt)$ in a similar way as in standard 
proofs of Fubini's theorem.)

\section{Invariant measures and disintegrations}\label{secterm}

Let $\bG$ be a locally compact 
second countable Hausdorff (in short {\em lcsc}) topological 
(multiplicative) group with unit element $e$.
The group $\bG$ is equipped with the Borel $\sigma$-field $\cG$.
Elements of $\bG$ will usually be denoted by $g$ or $h$. 
We fix a left-invariant locally finite {\em Haar measure}
$\lambda$ on $\bG$, see chapter 2 of \cite{Kallenberg} for more details 
and information. Left-invariance means
\begin{align*}
\int f(hg)\lambda(dg)=\int f(g)\lambda(dg),\quad h\in\bG, f\in\cG_+.
\end{align*}
The {\em modular function} is a continuous homomorphism 
$\Delta:\bG\rightarrow(0,\infty)$ satisfying
\begin{align}\label{modular}
\int f(gh)\lambda(dg)=\Delta(h^{-1})\int f(g)\lambda(dg),\quad h\in\bG,
\end{align}
for all $f\in\cS_+$. This modular function has the property
\begin{align}\label{modularinvers}
\int f(g^{-1})\lambda(dg)=\int \Delta(g^{-1})f(g)\lambda(dg),\quad f\in\cG_+.
\end{align}
The group $\bG$ is called {\em unimodular} if $\Delta(g)=1$
for all $g\in\bG$. By \eqref{modular} $\bG$ is unimodular 
if and only if $\lambda$ is right-invariant.

Let $(\bS,\cS)$ be a Borel space, i.e. a space Borel isomorphic to a Borel subset
of the interval $[0,1]$. 
Elements of $\bS$ will be named $s$ or $t$.
We assume that $\bG$ operates on $\bS$, i.e.\ we 
assume that there is a mapping $(g,s)\mapsto gs$
from $\bG\times\bS$ to $\bS$ having $g(hs)=(gh)s$ and $es=s$
for all $g,h\in\bG$ and $s\in\bG$. Here $e$ denotes the neutral
element of $\bG$.
The projections $\pi_s:\bG\rightarrow\bS$, $s\in\bS$, 
and the translations $\theta_g:\bS\rightarrow\bS$, $g\in\bG$, 
are given by 
$$
\pi_s(g)=\theta_g(s) = gs,\quad g\in\bG, s\in\bS.
$$ 
The set $\pi_s(\bG)=\bG s$ is called the {\em orbit} of $s$. 
We assume that the operation of $\bG$ 
on $\bS$ is {\em (measurably) proper}
in the sense that it is measurable as a map 
$\bG\times\bS\rightarrow \bS$ and that the set of all 
pushforwards $\mu_s:=\lambda\circ\pi_s^{-1}$, $s\in\bS,$ of the 
Haar measure under the projections is {\em uniformly} $\sigma$-finite. 
This means we require  the existence of a measurable 
partition $B_1,B_2,...$ of $\bS$
such that $\mu_s(B_n)<\infty$, $s\in\bS$, $n\in\N$.
This concept was introduced by Kallenberg in \cite{Kall07} 
and clearly generalizes the classical notion of a {\em topologically proper} 
operation of a lcsc group on a lcsc space (which is continous and where
$\pi_s^{-1}(K)$ is compact for any compact $K\subset S$ and all $s\in\bS$).
He also showed (\cite{Kall07}, Lemma 2.1) that properness is 
equivalent to the existence of 
a measurable function $k:\bS\rightarrow (0,\infty)$ such that 
\begin{align}\label{KallFunk}
\mu_sk = \int k(t)\mu_s(dt)<\infty,\quad s\in\bS.
\end{align}
Denote the cosets of the stabilizers as 
$$
G_{s,t}:=\{g\in\bG:gs=t\}=\pi_s^{-1}(\{t\}),\quad s,t\in\bS,
$$
which are measurable sets in $\bG$ under our assumptions on $\bS$. A measure $\nu$ 
on $\bS$ is called {\em invariant} (or $\bG$-{\em invariant}) if 
$$
\nu\circ\theta_g = \nu, \quad g\in\bG.
$$
The projection measures $\mu_s$, $s\in\bS$, are clearly 
invariant measures on $\bS$. They have the additional property that 
\begin{align}\label{mutrafo}
\mu_{gs} = \Delta(g^{-1})\mu_s,\quad g\in\bG, s\in\bS,
\end{align}
which means that the properness condition upon the operation enforces 
\begin{align}\label{deltaeins}
\Delta(g)=1,\quad g\in\bG_{s,s}, s\in\bS.
\end{align}
Hence the measures
\begin{align*}
\varphi_s:=\frac{\mu_s}{\mu_sk},\quad s\in\bS,
\end{align*}
are invariant, uniformly normalized in the sense that 
$\varphi_sk=1$, $s\in\bS$, 
and even constant on orbits, i.e. 
\begin{align}
\label{phiinv}\varphi_{gs}=\varphi_s, \quad s\in\bS, g\in\bG.
\end{align} 
By Fubinis theorem $\varphi$ is a kernel on $\bS$. 
Kallenberg proved in \cite{Kall07} that this kernel can be used as a normalized 
extremal generator of the 
convex cone of all $\sigma$-finite invariant measures on 
$\bS$ since for any such measure $\nu$ on 
$\bS$ 
\begin{align}\label{KallRep}
\nu(\cdot) = \int \varphi_s(\cdot)k(s)\nu(ds),
\end{align}
cf.\ Theorem 2.4 in \cite{Kall07}.

In the following Theorem \ref{lambdakern} we introduce a kernel $\kappa$ from $S\times S$ to $G$
that enables us to handle stabilizers and 
their cosets in $\bG$ within integral equations.
This kernel satisfies
\begin{align}
\label{explDesint}
\int f(gs,g)\lambda(dg) = \iint f(t,g)\kappa_{s,t}(dg)\mu_s(dt),
\quad f\in(\cS\otimes\cG)_+,s\in\bS.
\end{align}
In particular $\kappa$ disintegrates the Haar measure $\lambda$ on $\bG$  
along each orbit via
\begin{align}\label{explDesintHaar}
\int f(g)\lambda(dg) = \iint f(g)\kappa_{s,t}(dg)\mu_s(dt),
\quad f\in\cG_+,s\in\bS.
\end{align}

\begin{theorem}\label{lambdakern}
If $\bG$ operates properly on $\bS$ there is a 
kernel $\kappa$ from $S\times S$ to $\bG$ satisfying (\ref{explDesint}) 
and the following properties:
\begin{enumerate}
\item[{\rm (i)}] 
$\kappa_{s,gt}=\kappa_{s,t}\circ\theta_g^{-1},\quad g\in\bG,s,t\in\bS,$
\item[{\rm (ii)}] 
$\kappa_{s,t}$ is concentrated on $\bG_{s,t}:=\{g\in\bG:gs=t\}$ 
for $t\in\bG s,s\in\bS$,
\item[{\rm (iii)}] $\kappa_{s,t}(G)=1,\quad t\in\bG s, s\in\bS.$
\end{enumerate}
\end{theorem}

A kernel $\kappa$ with the above properties will be
fixed throughout the paper. In order to prove Theorem \ref{lambdakern} we will
need some more terminology and tools. When $\bG$ operates measurably on $\bS$ and $T$ 
we call a measure on a product 
space $S\times T$ {\em jointly $\bG$-invariant} if it is invariant 
with respect to the diagonal operation 
$$
\theta_g(s,t):=(gs,gt),\quad g\in\bG,s\in S,t\in T.
$$

Further, we call a kernel $\kappa$ from $S$ to $T$ {\em measurably} $\sigma${\em -finite} if for 
each $s\in\bS$ there is a measurable partition $B_1^s,B_2^s,...$ 
of $T$ such that $(s,t)\mapsto \I_{B_i^s}(t)$ is measurable for all 
$i\in\N$ and $\kappa(s,B_i^s)<\infty, s\in\bS$. It is easy to prove 
that a kernel from $S$ to $T$ is measurably $\sigma$-finite if and only 
if there exists a measurable function $f>0$ on $S\times T$ such 
that $\kappa_sf_s<\infty, s\in\bS,$ where $f_s:=f(s,\cdot)$.

Our aim is to disintegrate measurably 
labeled families of jointly $\bG$-invariant measures 
$\{M_r\}_{r\in R}$ on a product space in an invariant and 
measurable way. For this we need the following lemma which is a crucial extension
of well known results on the existence of
disintegrations of measures on product spaces (see e.g. \cite{Kallenberg} Theorem 6.3)
and their respective $\bG$-invariant versions for jointly $\bG$-invariant measures 
found by Kallenberg in \cite{Kall07}. Its proof is a straightforward adaption of
arguments found in \cite{Kallenberg} page 107 and \cite{Kall07} Theorem 3.5.

\begin{lemma}\label{messbareInvDesint}
Let $R,S,T$ be measurable spaces where $S$ and $T$ are Borel, 
$M$ be a measurably $\sigma$-finite kernel from $R$ to $S\times T$ 
and let $\bG$ operate measurably on both $S$ and $T$. 
\begin{enumerate}
\item[\rm (i)] 
There is a stochastic kernel $\nu$ from $R$ to $S$ and 
a measurably $\sigma$-finite kernel $\kappa$ from $R\times S$ to $T$ such that 
$$
M_r=\nu_r\otimes\kappa_r,\quad r\in R.
$$
\item[\rm (ii)] 
If $\nu'$ is a measurably $\sigma$-finite kernel from $R$ to 
$S$ with $M_r(\cdot\times T)\ll\nu_r', r\in R,$ there is a 
kernel $\kappa '$ from $R\times S$ to $T$ such that 
$$M_r=\nu_r'\otimes\kappa_r',\quad r\in R.$$
\item[\rm (iii)] 
If $M$ is such that $M_r$ is jointly $\bG$-invariant 
for each $r\in R$ and $\nu$ is a measurably $\sigma$-finite kernel from $R$ to $S$ 
such that $\nu_r$ is a $\bG$-invariant measure on $S$ 
with $M_r(\cdot\times T)\ll\nu_r$ for  $r\in R$ then there is a kernel 
$\kappa$ from $R\times S$ to $T$ with the invariance property 
$$
\kappa_r(gs,A) = \kappa_r(s,\theta_g^{-1}A),
\quad A\in\mathcal{T},s\in S,g\in\bG,r\in R,
$$
such that 
$$
M_r=\nu_r\otimes\kappa_r,\quad r\in R.
$$
\end{enumerate}
\end{lemma}
{\em Proof.}
{(i)} We may assume that 
$M_r(S\times T)>0, r\in R$. Since $M$ is measurably $\sigma$-finite 
we can choose a measurable function $f>0$ on $R\times S\times T$ such that 
$M_rf_r=1$, $r\in R$, and define the stochastic kernel $P$ from $R$ to 
$S\times T$ as $P_r:=f_r\cdot M_r$, $r\in R$. Then Proposition 7.26
in \cite{Kallenberg} yields a stochastic kernel $\tilde\kappa$ from $R\times S$ to $T$
such that together with the stochastic kernel $\nu_r:=P_r(\cdot\times T)$ 
$$
P_r=\nu_r\otimes\tilde\kappa_r,\quad r\in R,
$$ 
c.f. Dellacherie/Meyer \cite{DM} \RM{5}58. 
This is clearly equivalent to
$$
M_r=\nu_r\otimes\kappa_r,\quad r\in R,
$$ where $\kappa(r,s,A):=\int\I_A(t)/f(r,s,t)\tilde\kappa(r,s,dt)$, 
$A\in\mathcal{T}$, and thus proves the first assertion.

{(ii)}
If $\nu '$ is a given kernel from $R$ to $S$ with the property 
$M_r(\cdot\times T)\ll\nu_r', r\in R$, then 
$\nu_r\sim M_r(\cdot\times T)\ll\nu_r'$, $r\in R,$ and 
by Dellacherie/Meyer \RM{5}58 we may choose a measurable 
function $f:R\times S\rightarrow[0,\infty]$ such that 
$$
f(r,s) = \frac{d\nu_r}{d\nu_r'}(s), \quad\nu_r'\mbox{-a.e. } s\in S.
$$ 
Then 
$$
M_r=\nu_r'\otimes\kappa_r',\quad r\in R,
$$ 
where $\kappa '(r,s,\cdot):=f(r,s)\kappa(r,s,\cdot)$, which 
proves the second statement.

{(iii)}
From (ii) we get a kernel $\kappa$ from $R\times S$ to 
$T$ with 
$M_r = \nu_r\otimes\kappa_r, r\in R$. Invariance of $M_r$ and $\nu_r$ imply
for any $f\in(\cS\otimes\mathcal{T})_+$ that
\begin{align*}
\iint f(s,t)\kappa_{r,gs}(dt)\nu_r(ds) = 
\iint f(s,t)\kappa_{r,s}\circ\theta_g^{-1}(dt)\nu_r(ds),\quad 
g\in\bG,r\in R.
\end{align*}
Since $T$ is Borel this gives
$$
\kappa_{r,gs}=\kappa_{r,s}\circ\theta_g^{-1}, 
\quad\nu_r\mbox{-a.e. }s\in\bS,g\in\bG, r\in R.
$$ 
Fix some right Haar measure $\tilde\lambda$ on $\bG$. 
Fubinis theorem yields in particular
\begin{align}
\label{t1}
\kappa_{r,gs}=\kappa_{r,s}\circ\theta_g^{-1},
\quad \tilde\lambda\mbox{-a.e. }g\in\bG,\nu_r\mbox{-a.e. }s\in\bS,r\in R.
\end{align}
Let $l\geq 0$ be some measurable function on 
$\bG$ with $\tilde\lambda l=1$ and set
\begin{align*}
\overline\kappa_{r,s}:=\int(\kappa_{r,hs}\circ\theta_h)(l\cdot\tilde\lambda)(dh).
\end{align*}
A similar calculation as in \cite{Kall07} Theorem 3.5 shows that
on the sets
$$
A_r:=\{s\in\bS: \kappa_{r,ps}\circ\theta_p = 
\kappa_{r,qs}\circ\theta_q, \; \tilde\lambda^2\mbox{-a.e.}(p,q)\in\bG^2\},
\quad r\in R,
$$
we have 
\begin{align}\label{tmp1}
\overline\kappa_{r,s} = \overline\kappa_{r,hs}\circ\theta_h,\quad h\in\bG, s\in A_r, r\in R.
\end{align}
One now easily verifies by arguments relying on Fubinis theorem and a 
countable generator of $\cS$ that $(r,s)\mapsto \I_{A_r}(s)$ is measurable.
Further one can check that $A_r$ is $\bG$-invariant and \eqref{t1} implies that 
$\nu_r(A_r^c)=0$. Finally define
$$
\overline\kappa_{r,s}':=\I_{A_r}(s)\overline\kappa_{r,s},\quad s\in\bS,r\in R.
$$ 
Then by invariance of $A_r$ and \eqref{tmp1}
$$
\overline\kappa_{r,gs}'(A) =\overline\kappa_{r,s}'(\theta_g^{-1}A), 
\quad g\in\bG,s\in\bS,A\in\mathcal{T},r\in R,
$$
and since $\overline\kappa_{r,s}' = \overline\kappa_{r,s} = \kappa_{r,s},
\nu_r$-a.e. $s\in\bS$, the required disintegrations
$$
M_r = \nu_r\otimes\overline\kappa_r',\quad r\in R,
$$ 
hold indeed.
\qed

\vspace{0.3cm}

We are now ready to proof Theorem \ref{lambdakern}.

{\em Proof of Theorem \ref{lambdakern}.}
Consider the kernel 
$$
M_s:=\int \I\{(gs,g)\in\cdot\}\lambda(dg),\quad s\in\bS,
$$ 
from $S$ to $S\times\bG$ which is clearly measurably $\sigma$-finite by 
properness and has the property that every $M_s$ is a jointly $\bG$-invariant 
measure on $S\times\bG$. Further it is clear that 
$\mu_s:= \lambda\circ\pi_s^{-1}=M_s(\cdot\times G)$ 
and since the $\mu_s$ are $\sigma$-finite $\bG$-invariant measures we 
may apply Lemma \ref{messbareInvDesint} with $R:=S, T:=G$ and $\nu_s:=\mu_s$ 
to the kernel $M$ to obtain a kernel $\kappa$ from $\bS\times\bS$ to $\bG$ such that 
(\ref{explDesint}) and the invariance property (i) are fulfilled. 
It remains to show that $\kappa$ fulfills (ii),(iii): 
For (ii) note that for $s\in\bS$ by (\ref{explDesint})
\begin{align*}
\iint \I\{gs\not=t\}\kappa_{s,t}(dg)\mu_s(dt) 
&= \int \I\{gs\not=gs\}\lambda(dg) = 0.
\end{align*}
This means that $$\kappa_{s,t}(\bG_{s,t}^c)=0, \quad\mu_s\mbox{-a.e. } t\in\bS,s\in\bS,$$ 
and since $\mu_s\not=0$ for each $s\in\bS$ we may 
pick some $t\in\bG s$ such that $\kappa_{s,t}(\bG_{s,t}^c)=0$ holds. 
But then by (i) $\kappa_{s,t}(\bG_{s,t}^c)=0$ for all $t\in\bG s$. 
For (iii) choose $k$ as in \eqref{KallFunk} and note that setting $f(t,g):=k(t)$ in 
(\ref{explDesint}) yields 
\begin{align}
\mu_sk = \int k(t)\kappa_{s,t}(G)\mu_s(dt) = \kappa_{s,s}(G)\mu_sk,
\quad s\in\bS,
\end{align}
where we applied (i) in the last step and which 
implies $\kappa_{s,s}(G)=1=\kappa_{s,t}(G)$ for $t\in\bG s$ again by (i).
{\qed}

\vspace{0.3cm}

\begin{example}\label{exhomo}\rm Assume that $\bG$ operates
transitively on $\bS$, i.e.\ that there is only
one orbit. Fix some $c\in\bS$.
By \eqref{mutrafo} the measures $\mu_s$, $s\in\bS$,
are all multiples of $\mu_c$. By Corollary 2.6 in \cite{Kall07}
$\mu_c$ is up to normalization the unique invariant
$\sigma$-finite measure on $\bS$. The kernel $\kappa$ can be constructed 
by a suitable translation of the
probability measure $\kappa_c:=\kappa_{c,c}$.
Indeed, let $s\in\bS$ and take some $g_s\in\bG_{c,s}$, i.e.\ $g_sc=s$.
Then $G_{s,s}=g_sG_{c,c}g_s^{-1}$ and it is easy to see
that $\int\I\{g_sgg^{-1}_s\in\cdot\}\kappa_c(dg)$
is a left-invariant measure on $G_{s,s}$. If $G_{s,s}$ is compact,
then this measure must coincide with $\kappa_{s,s}$, see also
Corollary \ref{c32}. Now take $t\in\bS$ and some 
$g\in G_{s,t}$. By Theorem \ref{lambdakern} (ii) we then have
$\kappa_{s,t}=\kappa_{s,s}\circ\theta_g^{-1}$.
\end{example}

\begin{example}\label{lastex0}\rm We may further specialize
Example \ref{exhomo} by 
assuming that $\bS=\bG$ and that $(g,s)\mapsto gs$ is just the multiplication
in the group. Then $\mu_s=\Delta(s^{-1})\lambda$ for all $s\in\bS$.
For $s,t\in G$ we have $G_{s,t}=\{ts^{-1}\}$, while $\kappa_{s,t}$
is the Dirac measure located at $ts^{-1}$.
\end{example}

In applications, if some given operation is proper, it is 
usually not hard to determine a suitable partition or simultanously 
$\mu_s$-integrable function $k>0$ on $\bS$ and hence to actually prove properness. 
Conversely if this fails it can be hard to prove that a given 
operation is not proper.
The kernel $\kappa$ now gives a tool that enables us 
to reject properness in certain cases. Say that a subset $L\subset\bG$ is 
{\em locally closed} if it is the intersection of an open and a closed set. 
It is well known that such sets
inherit local-compactness from $\bG$, see also \cite{Bourbaki} I.65.

\begin{corollary}\label{c32}
Let $\bG$ operate properly on the Borel space $\bS$ such that 
$\bG_{s,s}$ is locally closed in $\bG$ for all $s\in\bS$. 
Then $\bG_{s,s}$ is compact in $\bG$ for all $s\in\bS$.
\end{corollary} 
{\em Proof.}
The assumption that $\bG_{s,s}$ is locally closed implies that 
$\bG_{s,s}$ is a locally compact subgroup of $\bG$ and for each $s$ we may choose
some left Haar measure $\lambda_s$ on $\bG_{s,s}$.
Consider the kernel $\kappa$ from Theorem $\ref{lambdakern}$. 
$\kappa_{s,s}$ is concentrated on $\bG_{s,s}$ and 
for any $g\in\bG_{s,s}$ we have by invariance
\begin{align*}
\kappa_{s,s}\circ\theta_g^{-1} = \kappa_{s,gs} = \kappa_{s,s}.
\end{align*}
The uniqueness theorem in \cite{Kall07} (Corollary 2.6) now implies 
$\lambda_s = c\cdot\kappa_{s,s}$, hence  
$\lambda_s$ is finite which implies compactness of $\bG_{s,s}$ 
(see \cite{Folland}, Proposition 11.4 d).\qed

\begin{example}\label{Beispieleproper}\rm Trivially the compact group $SO(d)$ operates
properly on $\R^d$ or on the affine Grassmanian $A(k,d)$ of $k$-dimensional affine subspaces 
of $\R^d$ (in fact on any space), and one readily proves that any $\R^u$ operates 
properly on $\R^d$ for $1\leq u \leq d$ via translation. On the other hand
by means of Corollary \ref{c32} 
it is straightforward to see that $\R^d$ does not operate properly on the 
Grassmanian $A(k,d)$ via translation. 
\end{example}

Choosing a system $\cO$ of representatives 
of the orbits $\bG s$, $s\in\bS$, 
the space $\bS$ splits into the disjoint union 
$$\bS = \bigcup_{b\in\cO} \bG b,$$ 
and we may consider the choice function 
$\beta:\bS\rightarrow\bS$ defined by $\beta(s):=b_s$, 
$b_s$ denoting the previously chosen representative of $\bG s$. A fixed 
choice of a system of representatives allows for the following 
canonical transfer of the modular function from $\bG$ to $\bS$ 
\begin{align}\label{lastDelta*}
\Delta^*(s):=\Delta(g^{-1}_s),\quad s\in\bS,
\end{align} 
where $g_s\in\bG_{\beta(s),s}$. This definition is independent of the
choice of $g_s$. Indeed, if $g,h\in\bG_{\beta(s),s}$ then
$g^{-1}h\in\bG_{\beta(s),\beta(s)}$ 
so that \eqref{deltaeins}  implies
$1=\Delta(g^{-1}h)$, i.e.\ $\Delta(g^{-1})=\Delta(h^{-1})$. 

As seen in Corollary \ref{c32} properness imposes restrictions upon the 
size of the stabilizers. This affects the relative size of the orbits. 
Most accessible is the case of countable $\bS$ which is of independent 
interest for applications (e.g.\ for percolation on countable graphs,
see \cite{B:L:P:S:99},\cite{LP:book}). 
The cardinality of a set $A$ is denoted by $|A|$.

\begin{lemma}\label{Scountable}
$\bG$ operates properly on a countable set $\bS$ if and only if
\begin{align*}
0<\lambda(G_{s,s})<\infty,\quad s\in\bS.
\end{align*}
In this case
\begin{align}\label{quothaar}
\Delta^*(s) = \frac{\lambda(G_{s,s})}{\lambda(G_{\beta(s),\beta(s)})}
= \frac{|G_{s,s}\beta(s)|}{|G_{\beta(s),\beta(s)}s|} ,\quad s\in\bS,
\end{align}
and either all orbits are infinite or all orbits are finite.
\end{lemma}
{\em Proof.}
Countability implies $0<\lambda(G_{s,s})$, $s\in\bS$, 
since $\lambda(G_{s,s})=0$ for some $s$
implies $\lambda(G) = \sum_{t\in\bG s}\lambda(\bG_{s,t}) = 0$ 
by left-invariance of $\lambda$ - an impossibility. For any $k$ on $\bS$
we have
\begin{align}\label{intmu}
\mu_sk = \int k(gs)\sum_{t\in\bG s}\I\{gs=t\}\lambda(dg) = 
\lambda(G_{s,s})\sum_{t\in\bG s}k(t),\quad s\in\bS.
\end{align}
Hence if the operation is proper then 
$\lambda(\bG_{s,s})<\infty$ for any $s$. This equation
also shows the converse since we may always choose $k>0$ on $\bS$ such that 
$\sum_{t\in\bG s}k(t)<\infty$, $s\in\bS$.
Now choose $k>0$ on $\bS$ such that 
$\sum_{s\in Gb}k(s)<\infty$, $b\in O$, and even 
$\sum_{s\in Gb}k(s)=\sum_{s\in Gb'}k(s)$, $b,b'\in O$. Then
\begin{align*}
\Delta^*(s) = \frac{\mu_sk}{\mu_{\beta(s)}k} = 
\frac{\lambda(G_{s,s})\sum_{t\in\bG s}k(t)}{\lambda(G_{\beta(s),\beta(s)})
\sum_{t\in\bG \beta(s)}k(t)} = 
\frac{\lambda(G_{s,s})}{\lambda(G_{\beta(s),\beta(s)})},\quad s\in\bS.
\end{align*}
By left-invariance of Haar measure we have 
\begin{align*}
\frac{\lambda(\bG_{s,s})}{\lambda(\bG_{s,s}\cap \bG_{t,t})} = 
[\bG_{s,s}:\bG_{s,s}\cap\bG_{t,t}] = |G_{s,s}t|,\quad s,t\in\bS,
\end{align*}
where $[\bG:H]$ denotes the left index of a subgroup $H\subset\bG$,
i.e.\ the number of the left cosets of $H$.
This yields the second identity in \eqref{quothaar}.
For any orbit $\bG s$ we have $\lambda(\bG)=|\bG s|\lambda(\bG_{s,s})$ 
and hence if 
$|\bG t|=\infty$ for some $t\in\bS$, then necessarily $\lambda(G)=\infty$ 
and thus for any other orbit $\bG s$ also $|\bG s|=\infty$ by properness.
\qed

\vspace{0.3cm}

\begin{remark}
Our hybrid setting of a lcsc group acting properly in a
measurable sense is taken from Kallenberg \cite{Kall07}.
It is more general than the usual assumption
of a topologically proper group operation and a 
first important step towards completely abandoning
topological assumptions both on $\bS$ and $\bG$.
However, we
have no substantial example falling in the more general
but not in the more specific category.
\end{remark}

If $\bS$ is not countable we need to establish measurability of $\beta$ and 
$\Delta^*$. For this recall the concept of universal measurability.
If $\mu$ is a measure on $(S,\cS)$ then $\cS^\mu$ denotes the completion of 
$\cS$ with respect to $\mu$.
The {\em universal completion} of a $\sigma$-algebra $\mathcal{S}$ is 
then defined as
\begin{align*}
\cS^u = \bigcap_\mu\cS^\mu
\end{align*}
where the intersection is taken over the class of all finite measures 
(or simply over the class of probability measures)
on $(\bS,\cS)$. Its elements are called {\em universally measurable} 
sets and a map $f:\bS\rightarrow \bT$ is called 
{\em universally measurable} if it is 
$\cS^u/\cT$-measurable.

\begin{lemma}
\label{univmessbar}
The following holds:
\begin{enumerate}
\item[{\rm (i)}] The orbits 
$\bG s$, $s\in\bS,$ are universally measurable sets in $\bS$;
\item[{\rm (ii)}] If $O\in\cS$ then $\beta$ is universally measurable;
\item[{\rm (iii)}]
If $O\in\cS$ then for any $B\in\cG$ the map $s\mapsto\kappa_{\beta(s),s}(B)$ 
is universally measurable;
\item[{\rm (iv)}] If $O\in\cS$ then $\Delta^*$ is universally measurable. 
\end{enumerate} 
\end{lemma}
{\em Proof.} Consider the Borel isomorphism 
$\psi:\bG\times\bS\rightarrow\bG\times\bS$ 
given by $\psi(g,s)=(g,gs)$ and the measurable sets 
$A_s:=\psi(\bG\times\{s\})$, $s\in\bS$. 
Since $\bG$ is Borel the projection of $A_s$ on $\bS$ 
is a universally measurable set in $\bS$ according to Dellacherie and 
Meyer \cite{DM} Section \RM{3}.44. These projections are clearly the 
orbits of the operation, hence (i) follows. 

For (ii) note that for $B\subset S$ we have $\beta^{-1}(B)=\bG (B\cap O)$. 
Hence $\beta^{-1}(B)= G(B\cap O) = \pr_S(\varphi(G\times (B\cap O)))$ and for 
$B,O\in\cS$ this implies that $\beta^{-1}(B)$ is universally measurable 
since $G$ is Borel.

Assertion (iii) holds since the map $(s,t)\mapsto\kappa_{s,t}(B)$ is measurable according to 
Theorem \ref{lambdakern} while $s\mapsto(\beta(s),s)$ is universally 
measurable by (ii) and elementary properties of 
the product $\sigma$-algebra.

The universal measurability of $\Delta^*$ and hence (iv) follows from the
representation 
\begin{align}\tag*{\qed}
\Delta^*(s) = \int \Delta(g^{-1})\kappa_{\beta(s),s}(dg),\quad s\in\bS.
\end{align}

\begin{remark}\rm
The concept of universal measurability is useful 
when dealing with finite or at least $\sigma$-finite measures $\mu$ on 
a measurable space $(\bS,\cS)$ since for a universally measurable map 
$f$ integrals $\mu f$ with respect to any such $\mu$ make 
sense. Any $\sigma$-finite $\mu$ has a unique extension to the 
class of $\cS^\mu$-measurable functions and in particular to the (smaller) 
class of $\cS^u$-measurable functions. We make heavy use of this fact for 
almost all results in this paper without further notice.
\end{remark}

Throughout this paper we fix one system $O$ of representatives of 
the orbits and require $O\in\cS$ such that $\beta$ and $\Delta^*$
are universally measurable by Lemma \ref{univmessbar}. 
By means of (\ref{phiinv}) and $\beta$ equation (\ref{KallRep}) 
can be modified as follows:
\begin{align}\label{OrbitRep}
\nu(\cdot) &= \int \varphi_s(\cdot)k(s)\nu(ds) = 
\int\mu_{\beta(s)}(\cdot)\mu_{\beta(s)}(k)^{-1}k(s)\nu(ds) \notag \\
&=\int \mu_b(\cdot)\nu^*(db),
\end{align}
where $\nu^*= (\mu_{\beta}(k)^{-1}k\cdot\nu)\circ\beta^{-1}$ 
is a measure concentrated on 
$\cO$, in the sense that any measurable $B\subset\bS$ being 
disjoint with $\cO$ has $\nu^*(B)=0$.

\section{Random measures and Palm pairs}\label{secPalm}

As before let $(\bS,\cS)$ be a Borel space.
Let $\bM(\bS)$ denote the space of all $\sigma$-finite  measures on $\bS$.
We equip $\bM(\bS)$ with the smallest $\sigma$-field
$\cM(\bS)$ making the mappings $\mu\mapsto \mu(B)$ for all
$A\in\cS$ measurable.
Let $(\Omega,\cA,\BP)$ be a $\sigma$-finite measure space. We
use a probabilistic language even though $\BP$ need not be
a probability measure. A {\em random measure} on $\bS$ is a
measurable mapping $\xi:\Omega\rightarrow \bM(\bS)$ that
is uniformly $\sigma$-finite in the sense that there is a
partition $B_1,B_2,...$ of $\bS$ such that
$\xi(B_i)<\infty$ $\BP$-a.e.\ for any $i\in\N$.
We use the kernel notation $\xi(\omega,B):=\xi(\omega)(B)$.
If $\xi$ is a random measure on $\bS$ then the {\em Campbell measure} 
$C_\xi$ of $\xi$ (with respect to $\BP$) is the measure on $\Omega\times\bS$ 
satisfying
\begin{align*}
C_\xi f = \iint f(\omega,s)\xi(\omega,ds)\BP(d\omega),\quad 
f\in(\cA\otimes\cS)_+.
\end{align*}
This measure is finite on sets of the form
$\{\omega\in A:\xi(\omega,B_i)\le n\}\times B_i$,
where $i,n\in \N$, $A\in\cA$ has $\BP(A)<\infty$ and
the $B_i$ are as in the definition of $\xi$.
It follows that $C_\xi$ is $\sigma$-finite. 
Hence there exists a {\em supporting measure} of $\xi$,
i.e.\ a $\sigma$-finite measure $\nu$ on
$\bS$ such that $C_\xi(\Omega\times \cdot)$ and $\nu$ are equivalent in the sense 
of mutual absolute continuity. 
If $(\Omega,\cA)$ is Borel, then there
is a $\sigma$-finite kernel $Q$ from $\bS$ to $\Omega$ 
disintegrating $C_\xi$ as follows:
\begin{align}\label{Palm}
C_\xi f = \iint f(\omega,s)Q_s(d\omega)\nu(ds),\quad f\in(\cA\otimes\cS)_+.
\end{align}
We call a pair $(\nu,Q)$ satisfying \eqref{Palm} a {\em Palm pair} of $\xi$. 
The kernel $Q$ is the $\nu$-associated {\em Palm kernel} 
of $\xi$. To make the dependence on $\xi$ (but not on $\BP$)
explicit, we often write $(\nu_\xi,Q_\xi):=(\nu,Q)$.

We now return to the setting established in Section \ref{secterm}.
In addition we assume that $\bG$ is operating measurably
on $\Omega$. There is no risk of confusion to denote
the associated translations by $\theta_g:\Omega\rightarrow\Omega$, $g\in\bG$. 
The set $\{\theta_g\}_{g\in\bG}$ is commonly 
refered to as a {\em flow} on $\Omega$. 
Our basic assumption is that $\BP$ is invariant under the flow,
i.e.
\begin{align}\label{stat}
\BP\circ\theta^{-1}_g =\BP, \quad g\in G.
\end{align}
A random measure $\xi$ on $\bS$  is called {\em invariant} 
(or {\em $\bG$-invariant}) if it satisfies 
\begin{align}\label{xicov}
\xi(\theta_g\omega,B) = \xi(\omega,g^{-1}B),
\quad g\in\bG,\omega\in\Omega,B\in\cS.
\end{align}
Similarly a kernel $Q$ from $\bS$ to $\Omega$ is called 
{\em invariant}  if
\begin{align}\label{Qcov}
Q (gs,A) = Q(s,\theta_g^{-1}A), \quad g\in G,s\in\bS, A\in\mathcal{A}.
\end{align}
Assume that $\xi$ is an invariant random measure.
Then it is easy to see that the Campbell measure of
$\xi$ is jointly invariant. 
If $(\Omega,\cA)$ is Borel, Corollary 3.5 in \cite{Kall07} implies 
that there is an {\em invariant Palm pair} $(\nu,Q)$ of $\xi$,
meaning that both $\nu$ and $Q$ are invariant.
Next we formulate the {\em refined Campbell theorem}. Although
a simple consequence of the definitions, it
is the main tool of Palm calculus for invariant random measures.
Recall the representation
\eqref{OrbitRep} and that $\theta_e$ is the identity on $\Omega$.
Adapting common terminology of probability theory,
we denote integration with respect to a measure
$Q'$ on $\Omega$ by $\BE_{Q'}$ and the $\theta_g$'s may
be interpreted as random elements of $\Omega$ in the following. 

\begin{proposition}\label{refC}
Assume that the invariant random measure $\xi$
has an invariant Palm pair $(\nu,Q)$. Then, 
for all $f\in(\cA\otimes{\cG}\otimes\cS)_+$,
\begin{align}\label{C}
\BE_\BP\iint f(\theta_g^{-1},g,\beta(t))\kappa_{\beta(t),t}(dg)
\xi(dt)
=\int \BE_{Q_b}\int f(\theta_e,g,b)\lambda(dg)\nu^*(db).
\end{align}
\end{proposition}
{\em Proof.} Let $f\in(\cA\otimes{\cG}\otimes\cS)_+$ and denote
the right-hand side of \eqref{C} by $I$.
By \eqref{explDesintHaar} and Fubini's theorem,
$$
I=\iint \BE_{Q_b}\int f(\theta_e,g,b)\kappa_{b,t}(dg)\mu_b(dt)\nu^*(db).
$$
By Theorem \ref{lambdakern} (ii),  $\kappa_{b,t}$ is concentrated
on $G_{b,t}$. Hence we obtain by invariance \eqref{Qcov} of the
Palm kernel that
$$
I=\iint \BE_{Q_t}\int f(\theta^{-1}_g,g,b)\kappa_{b,t}(dg)\mu_b(dt)\nu^*(db).
$$
For $b\in O=\beta(S)$ and $t\in Gb$ we have $b=\beta(t)$.
Since $\nu^*$ is concentrated on $O$ 
we obtain 
$$
I=\iint \BE_{Q_t}\int f(\theta^{-1}_g,g,\beta(t))
\kappa_{\beta(t),t}(dg)\mu_b(dt)\nu^*(db).
$$
An application of \eqref{OrbitRep} yields
$$
I=\int \BE_{Q_t}\int f(\theta^{-1}_g,g,\beta(t))
\kappa_{\beta(t),t}(dg)\nu(dt).
$$
The defining property \eqref{Palm} of a Palm pair yields
the assertion \eqref{C}.\qed

\begin{remark}\label{remgen}\rm
Assume that $\bG=\{e\}$. Then $\lambda=\delta_e$, $O=S$,
$\mu_s=\delta_s$ and $\kappa_{s,t}=\I\{s=t\}\delta_e$.
Let $\xi$ be a random measure on $\bS$ with supporting
measure $\nu$. Then $\xi$ is invariant. Let $(\nu,Q)$ be a Palm 
pair of $\xi$.  Then $(\nu,Q)$ is invariant, $\nu^*=\nu$, and  
the refined Campbell theorem \eqref{C} boils down
to the defining equation \eqref{Palm} of a Palm pair.
Hence general random measures can be treated within our
framework of invariant random measures. 
\end{remark}

\begin{example}\label{exhomo2}\rm 
Consider the situation of Example \ref{exhomo}, i.e.\ assume
that $\bG$ operates transitively on $\bS$. Let $\xi$ be
an invariant random measure on $\bS$. Fixing some $c\in\bS$, 
we can take $\nu:=\mu_c$ as a supporting measure
of $\xi$. Moreover, taking $\beta\equiv c$, we 
clearly have $\nu^*=\delta_c$. Then \eqref{C}
simplifies to 
\begin{align}\label{Chomo}
\BE_\BP\iint f(\theta_g^{-1},g)\kappa_{c,t}(dg)\xi(dt)
=\BE_{Q_c}\int f(\theta_e,g)\lambda(dg),\quad
f\in(\cA\otimes{\cG})_+.
\end{align}
In particular,
\begin{align}\label{Palmhomo}
Q_c(A)=\BE_\BP \iint \I_A(\theta_g^{-1})w(t)\kappa_{c,t}(dg)\xi(dt),
\quad A\in\cA,
\end{align}
where $w\in \cS_+$ has $\int wd\mu_c=1$.
In fact one can use \eqref{Palmhomo} to give an explicit
definition of $Q_c$ (without any Borel assumption on $\Omega$)
and then derive \eqref{Chomo}.
This is the approach taken in \cite{RoZae90} and \cite{La08b}.
In the further special case $\bS=\bG$ of Example \ref{lastex0}
we may take $c=e$ and \eqref{Chomo} simplifies to
\begin{align}\label{Cgroup}
\BE_\BP\int f(\theta_g^{-1},g)\xi(dg)
=\BE_{Q_e}\int f(\theta_e,g)\lambda(dg),\quad
f\in(\cA\otimes{\cG})_+.
\end{align}
Equation \eqref{Palmhomo} changes accordingly.
This {\em skew factorization} is the standard approach to
Palm measures of stationary random measures on a group.
We refer to \cite{Kall07} for more details and
historical remarks.
\end{example}

\begin{example}\label{exquasit}\rm 
Assume that $O$ is finite, i.e.\ that there are only finitely
many orbits. Following  \cite{B:L:P:S:99} we call this
the {\em quasi-transitive} case. Let $\xi$ be an invariant random measure on $\bS$.
Then without loss of generality a supporting measure may be chosen as
$\nu:=\sum_{b\in O'}\mu_b$ for some $O'\subset O$. 
It is easy to check that $\nu^*=\sum_{b\in O'}\delta_b$.
Then if Q is an associated invariant Palm kernel of $\xi$
the refined Campbell theorem \eqref{C} now implies
for all $f\in(\cA\otimes{\cG})_+$ that
\begin{align}\label{Cquasi}
\BE_\BP\iint f(\theta_g^{-1},g)\kappa_{b,t}(dg)\xi_b(dt)
=\BE_{Q_b}\int f(\theta_e,g)\lambda(dg),\quad b\in O',
\end{align}
where $\xi_b$ is the restriction of $\xi$ to the orbit $Gb$.
Let $b\in O'$ and $w_b\in\mathcal{S}_+$ with $\int w_bd\mu_b=1$.
Then \eqref{Cquasi} implies
\begin{align}\label{Palmquasi}
Q_b(A)=\BE_\BP \iint \I_A(\theta_g^{-1})w_b(t)\kappa_{b,t}(dg)\xi_b(dt),
\quad A\in\cA.
\end{align}
This can be used for defining the Palm kernel $Q$ explicitly, just as
in the transitive case of Example \ref{exhomo2}.
\end{example}


\section{The transport formula}\label{sexch}

In the remainder of the paper we assume that the lcsc group $G$ 
operates measurably on $(\Omega,\cA)$ and properly on 
the Borel space $(\bS,\cS)$.  
In this section we fix an invariant $\sigma$-finite measure $\BP$ on
$(\Omega,\cA)$. Our aim is to 
derive a fundamental transport property of Palm measures.
In the special case where $\bG=\bS$ is an Abelian group
the result boils down to Theorem 3.6 in \cite{LaTh08}.
Other special cases will be discussed below.
We use the function $\Delta^*$ defined by \eqref{lastDelta*}.
A kernel $\gamma$ from $\Omega\times\bS$ to $\bS$ is called {\em invariant} if
\begin{align}\label{equivar}
\gamma(\theta_g\omega,gs,B)=\gamma(\omega,s,g^{-1}B),
\quad g\in\bG,s\in\bS,\omega\in\Omega,B\in\cS.
\end{align}

\begin{theorem}\label{thlastcycle}  
Consider two invariant random measures 
$\xi$ and $\eta$ on $\bS$ and let $\gamma$ and $\delta$ be 
invariant kernels from $\Omega\times\bS$ to $\bS$ satisfying 
\begin{align}\label{lastT*}
\iint\I\{(s,t)\in\cdot\}\gamma(\omega,s,dt)\xi(\omega,ds)
=\iint\I\{(s,t)\in\cdot\}\delta(\omega,t,ds)\eta(\omega,dt)
\end{align}
for $\BP$-a.e.\ $\omega\in\Omega$ and $(\nu_\xi,Q_\xi)$, $(\nu_\eta,Q_\eta)$ be invariant Palm pairs
of $\xi$ and $\eta$ respectively. Then we have for any measurable function 
$f\in (\cA\otimes{\cG}\otimes\cS\otimes\cS)_+$ that
\begin{align}\label{Neveutrans}
\int\BE_{Q_{\eta,b}}&\iint f(\theta_g^{-1},g^{-1},b,\beta(s))
\Delta^*(s)\kappa_{\beta(s),s}(dg)\delta(b,ds)\nu^*_\eta(db) \notag \\
 &=\int\BE_{Q_{\xi,b}}\iint f(\theta_e,g,\beta(s),b)
\kappa_{\beta(s),s}(dg)\gamma(b,ds)\nu^*_\xi(db).
\end{align}
\end{theorem}
{\em Proof.}
Let $k>0$ be as in (\ref{KallFunk}). 
Then for any $b\in O=\beta(S)$ and $g\in\bG$
\begin{align*}
l(b)\int k(g^{-1}hb)\lambda(dh)=1,
\end{align*}
where $l(b):=\mu_b(k)^{-1}$. Take 
$f\in (\cA\otimes{\cG}\otimes\cS\otimes\cS)_+$
and denote the right-hand side of \eqref{Neveutrans} by $I$.
By Fubini's theorem,
\begin{align*}
I=\int\BE_{Q_{\xi,b}}\iiint f(\theta_e,g,\beta(s),b)l(b)
k(g^{-1}hb)\kappa_{\beta(s),s}(dg)\gamma(b,ds)\lambda(dh)\nu^*_\xi(db).
\end{align*}
Applying the refined Campbell theorem \eqref{C} gives that $I$ equals
\begin{align*}
\BE_{\BP}\iiiint & f(\theta^{-1}_h,g,\beta(s),\beta(t))l(\beta(t))
k(g^{-1}h\beta(t))\\
&\qquad\qquad\qquad\kappa_{\beta(s),s}(dg)
\gamma(\theta^{-1}_h,\beta(t),ds)\kappa_{\beta(t),t}(dh)\xi(dt)\\
&=\BE_{\BP}\iiiint f(\theta^{-1}_h,g,\beta(s),\beta(t))l(\beta(t))k(g^{-1}t)\\
&\qquad\qquad\qquad\kappa_{\beta(s),h^{-1}s}(dg)
\gamma(h\beta(t),ds)\kappa_{\beta(t),t}(dh)\xi(dt)\\
&=\BE_{\BP}\iiiint f(\theta^{-1}_h,g,\beta(s),\beta(t))l(\beta(t))k(g^{-1}t)\\
&\qquad\qquad\qquad\kappa_{\beta(s),h^{-1}s}(dg)
\kappa_{\beta(t),t}(dh)\gamma(t,ds)\xi(dt),
\end{align*}
where we have the invariance property \eqref{equivar} of $\gamma$, 
invariance of $\beta$
and the fact that $\kappa_{\beta(t),t}$ is concentrated on
$G_{\beta(t),t}$ (see Theorem \ref{lambdakern} (ii)).
By Theorem \ref{lambdakern} (i) and \eqref{lastT*} 
\begin{align*}
I=\BE_{\BP}\iiiint f(\theta^{-1}_h,h^{-1}g,\beta(s),\beta(t))&l(\beta(t))
k(g^{-1}ht) \\
&\kappa_{\beta(t),t}(dh)\delta(s,dt)\kappa_{\beta(s),s}(dg)\eta(ds).
\end{align*}
Using the invariance properties of the kernels $\delta$ and $\kappa$,
we obtain that $I$ equals
\begin{align*}
\BE_{\BP}\iiiint f(\theta^{-1}_h\circ\theta^{-1}_g,&h^{-1},\beta(s),\beta(t))
l(\beta(t))\\
&k(hgt)\kappa_{\beta(t),t}(dh)\delta(\theta^{-1}_g,\beta(s),dt)
\kappa_{\beta(s),s}(dg)\eta(ds),
\end{align*}
where we have used that $\theta_{gh}^{-1}=\theta^{-1}_h\circ\theta^{-1}_g$
and that $g^{-1}s=\beta(s)$ for $s,g$ as in the above integral.
At this stage we can use the refined Campbell theorem \eqref{C}
for $\eta$ to obtain that $I$ equals
\begin{align*}
\int\BE_{Q_{\eta,b}}\iiint f(\theta^{-1}_h,h^{-1},b,\beta(t))
l(\beta(t))k(hgt)\kappa_{\beta(t),t}(dh)\delta(b,dt)\lambda(dg)\nu^*_\eta(db).
\end{align*}
Now take $h\in\bG$ and $t\in\bS$ with $h\beta(t)=t$. Then
$$
\int k(hgt)\lambda(dg)=\int k(gh\beta(t))\lambda(dg)
=\Delta(h^{-1}) l(\beta(t))^{-1}.
$$
Hence we obain from Fubini's theorem that
$I$ equals the left-hand side of \eqref{Neveutrans}.\qed

An immediate consequence of Theorem \ref{thlastcycle} is the following 
{\em exchange formula}  for Palm pairs.
A first version of this fundamental and very useful
formula was obtained by Neveu (see e.g.\ \cite{Neveu}).

\begin{corollary}\label{NeveuFormula}
Let $\xi$ and $\eta$ be invariant random measures on $\bS$. 
Assume that $\xi$ and $\eta$ admit invariant Palm pairs
$(\nu_\xi,Q_\xi)$ and $(\nu_\eta,Q_\eta)$, respectively. 
Then for any $f\in (\cA\otimes{\cG}\otimes\cS\otimes\cS)_+$
\begin{align}\label{Neveu}
\int\BE_{Q_{\eta,b}}\iint &f(\theta_g^{-1},g^{-1},b,\beta(s))
\Delta^*(s)\kappa_{\beta(s),s}(dg)\xi(ds)\nu^*_\eta(db) \notag \\
 &=\int\BE_{Q_{\xi,b}}\iint f(\theta_e,g,\beta(s),b)
\kappa_{\beta(s),s}(dg)\eta(ds)\nu^*_\xi(db).
\end{align}
\end{corollary}

Before discussing some examples we 
mention one consequence of Corollary \ref{NeveuFormula}
arising for a special choice of $f$.

\begin{corollary}\label{NeveuFormula2}
Under the hypothesis of Corollary \ref{NeveuFormula} we have
\begin{align}\label{Neveu2}
\int\BE_{Q_{\eta,b}}\iint  &f(\theta_g^{-1},g^{-1}b,\beta(s))
\Delta^*(s)\kappa_{\beta(s),s}(dg)\xi(ds)\nu^*_\eta(db) \notag \\
 &=\int\BE_{Q_{\xi,b}}\int f(\theta_e,s,b)\eta(ds)\nu^*_\xi(db),
\quad f\in (\cA\otimes\cS\otimes\cS)_+.
\end{align}
\end{corollary}
{\em Proof.}
Take $f\in (\cA\otimes\cS\otimes\cS)_+$ and apply
Corollary \ref{NeveuFormula} with the function
$\tilde f(\omega,g,s,t):=f(\omega,gs,t)$.\qed

\begin{example}\label{exchangegroup}\rm
Assume that $\bS=\bG$ as in Example \ref{lastex0}
and let $\xi$ and $\eta$ be invariant random measures on $\bG$.
Then \eqref{Neveu} means for $f\in (\cA\otimes{\cG})_+$
\begin{align}\label{Neveugroup}
\BE_{Q_{\eta}}\int f(\theta_g^{-1},g^{-1})\Delta(g^{-1})\xi(dg)
=\BE_{Q_{\xi}}\int f(\theta_e,g)\eta(dg),
\end{align}
where $Q_{\xi}:=Q_{\xi,e}$, $Q_{\eta}:=Q_{\eta,e}$, cf.\ Example \ref{exhomo2}.
For an Abelian group this is Neveu's exchange formula, see \cite{Neveu}.
\end{example}

\begin{example}\label{exchangequasi}\rm
Consider the quasi-transitive case of  Example \ref{exquasit}
and let $\xi,\eta,\gamma,\delta$ be as in Theorem \ref{thlastcycle}.
Then \eqref{Neveutrans} easily implies for all $f\in (\cA\otimes{\cG})_+$
and all $b,b'\in O$ that
\begin{align}\label{exquasi}\notag
\nu^*_\eta\{b\}\BE_{Q_{\eta,b}}\iint &f(\theta_g^{-1},g^{-1})\Delta^*(s)
\kappa_{b',s}(dg)\I\{s\in Gb'\}\delta(b,ds)\\
&=\nu^*_\xi\{b'\}\BE_{Q_{\xi,b'}}
\iint f(\theta_e,g)\kappa_{b,s}(dg)\I\{s\in Gb\}\gamma(b',ds).
\end{align}
The transitive special case of this result can be found
in \cite{La08b}. In case $\bG=\bS$ is an Abelian group
we recover Theorem 3.6 in \cite{LaTh08}
(see also \cite{La08} for the non Abelian case).
\end{example}

\section{Characterization of Palm pairs}

We consider a kernel $\xi$ from $\Omega$ to $\bS$ that is
uniformly $\sigma$-finite in the sense, that there is a
partition $B_1,B_2,...$ of $\bS$ such that
$\xi(B_i)<\infty$ for any $i\in\N$.
In contrast to the previous section we do not fix the 
underlying measure $\BP$ on $(\Omega,\mathcal{A})$. 
Instead we fix a $\sigma$-finite measure 
$\nu$ on $\bS$ and a $\sigma$-finite kernel $Q$ from $\bS$ to $\Omega$
and ask for conditions that are necessary and sufficient for $(\nu,Q)$
to be a Palm pair of $\xi$ with respect to some
$\sigma$-finite measure $\BP$. A first result in this direction
can be formulated without any invariance assumptions
which can be seen as a special case (see Remark \ref{remgen}) 
of some independent interest.

\begin{proposition}\label{CharPalm}
The pair  $(\nu,Q)$ is a Palm pair of $\xi$ with respect to some 
$\sigma$-finite  measure on $(\Omega,\mathcal{A})$
iff $\nu\otimes Q$ is $\sigma$-finite, $Q_s(\xi=0)=0$ for $\nu$-a.e. 
$s\in S$ and 
\begin{align}\label{Sym}
\iiint f(\omega,s,t)\xi(\omega,dt)Q_s(d\omega)\nu(ds) 
= \iiint f(\omega,t,s)\xi(\omega,dt)Q_s(d\omega)\nu(ds)
\end{align}
for any $f\in(\mathcal{A}\otimes\mathcal{S}\otimes\mathcal{S})_+$. 
If $(\nu,Q)$ and $\xi$ are invariant the same characterisation holds 
and in addition the underlying measure on $(\Omega,\mathcal{A})$ may be chosen
to be invariant.
\end{proposition}
{\em Proof.} First, assume that $(\nu,Q)$ is a Palm pair of $\xi$ with respect 
to some $\sigma$-finite measure $\BP$ on $\Omega$. This means 
\begin{align}\label{Camp}
\nu\otimes Q = C_\xi
\end{align}
where $C_\xi$ is the Campbell measure of $\xi$ w.r.t.\ $\BP$.
As we have seen earlier, $C_\xi$ is $\sigma$-finite. Further
\begin{align*}
\int Q_s(\xi=0)\nu(ds) = 
\BE_\BP\int \I\{\xi=0\}\xi(ds)= 0,
\end{align*}
which implies $Q_s(\xi=0)=0$ for $\nu$-a.e.\ $s\in S$. 
To show $\eqref{Sym}$ we take
$f\in(\mathcal{A}\otimes\mathcal{S}\otimes\mathcal{S})_+$
and obtain from \eqref{Camp} and Fubini's theorem
that the left hand-side of $\eqref{Sym}$ equals
\begin{align*} 
\BE_\BP\iint f(\theta_e,s,t)\xi(dt)\xi(ds)
&= \BE_\BP\iint f(\theta_e,s,t)\xi(ds)\xi(dt).
\end{align*}
Applying  \eqref{Camp} again, we see that the latter
expression coincides with the right-hand side of \eqref{Sym}.

We now prove the converse implication. By $\sigma$-finiteness 
we may choose a measurable function $g'>0$ on $\Omega\times \bS$ 
such that $(\nu\otimes Q)g'<\infty$. 
Since $\xi$ is uniformly $\sigma$-finite, we may choose 
$\tilde g>0$ on $\Omega\times \bS$ with 
$0<\xi(\tilde g)<\infty$ on $\{\xi\not=0\}$. 
Now set $g:=g'\wedge\tilde g$ and
$h(\omega,s):=g(\omega,s)/(\xi g)(\omega)$, 
where $h(\omega,s):=0$ if $\xi(g)=0$.
Define the measure $\BP$ by
\begin{align}\label{DefP}
\BP(A):=\iint \I_A(\omega)h(\omega,s)Q_s(d\omega)\nu(ds),
\quad A\in\mathcal{A}.
\end{align}
By assumption on $Q$ we have $\BP(\xi=0)=0$.  
The function  $\xi(g)$ is finite and positive on $\{\xi\not=0\}$.
Furthermore,
\begin{align*}
\BE_\BP\xi(g) = \iint (\xi g)(\omega)\frac{g(\omega,s)}{(\xi g)(\omega)} 
Q_s(d\omega)\nu(ds) \leq \iint g'(\omega,s)Q_s(d\omega)\nu(ds)<\infty.
\end{align*}
Hence $\BP$ is $\sigma$-finite. Moreover we have for
$f\in(\mathcal{A}\otimes\mathcal{S})_+$
\begin{align*}
\BE_\BP \int f(\theta_e,t)\xi(dt)
&=\iiint f(\omega,t)h(\omega,s)\xi(\omega,dt)Q_s(d\omega)\nu(ds)\\
&=\iiint f(\omega,s)h(\omega,t)\xi(\omega,dt)Q_s(d\omega)\nu(ds),
\end{align*}
where we have used $\eqref{Sym}$ to get the second identity. 
This is just $(\nu\otimes Q) f$ since $\xi(h)=1$ on
$\{\xi\not=0\}$ by definition of $h$, and 
$Q_s(\xi=0)=0$ for $\nu$-a.e.\ $s\in S$.

It remains to show that the measure $\BP$ defined in \eqref{DefP} 
is invariant for invariant $\nu$, $Q$, and $\xi$. Take
$f\in\mathcal{A}_+$ and $g\in\bG$. By invariance of $Q$
and $\nu$,
\begin{align*}
\BE_\BP f\circ\theta_g &=
\iint f(\theta_g\omega)h(\omega,s)Q_s(d\omega)\nu(ds) \\
&=\iint f(\omega)h(\theta_g^{-1}\omega,s)Q_{gs}(d\omega)\nu(ds) \\
&=\iint f(\omega)h(\theta_g^{-1}\omega,g^{-1}s)Q_{s}(d\omega)\nu(ds).
\end{align*}
Since $(\nu,Q)$ is a Palm pair of the invariant $\xi$ (w.r.t.\ $\BP$)
we obtain 
\begin{align*}
\BE_\BP f\circ\theta_g 
&=\BE_\BP \int f(\theta_e)h(\theta_g^{-1},g^{-1}s)\xi(ds)\\
&=\iint f(\omega)h(\theta_g^{-1}\omega, s)\xi(\theta_g^{-1}\omega, ds)
\BP(d\omega)=\BE_\BP f,
\end{align*}
where we have used in the last step that 
$\int h(\theta_g^{-1}\omega,s)\xi(\theta_g^{-1}\omega,ds)=1$ 
for $\BP$-a.e. $\omega$ since $\{\xi\not=0\}$ is invariant and has 
a complement of $\BP$-measure $0$.\qed

\vspace{0.3cm}

The  following main result of this section is
a significant extension of Mecke's \cite{Mecke} famous
characterization of Palm measures of stationary random
measures on an Abelian group. In the transitive special case
of Example \ref{exhomo2} the result has been derived in
\cite{RoZae90} and \cite{La08b}.

\begin{theorem}\label{charMecke}
Assume that $\xi$ and  $(\nu,Q)$ are invariant. 
Then $(\nu, Q)$ is a Palm pair of $\xi$ with 
respect to some invariant $\sigma$-finite measure on 
$(\Omega,\mathcal{A})$ 
iff $\nu\otimes Q$ is $\sigma$-finite, $Q_s(\xi=0)=0$ $\nu$-a.e.\
$s\in\bS$, and, 
for any $f\in(\mathcal{A}\otimes\mathcal{S}\otimes\mathcal{S})_+$,
\begin{multline}\label{MeckeGl}
\int\BE_{Q_b}\iint f(\theta_g^{-1},g^{-1}b,\beta(s))\Delta^*(s)
\kappa_{\beta(s),s}(dg)\xi(ds)\nu^*(db)\\
=\int\BE_{Q_b}\int f(\theta_e,s,b)\xi(ds)\nu^*(db).
\end{multline}
\end{theorem}

{\em Proof.}
If $(\nu,Q)$ is a Palm pair of $\xi$ then $\sigma$-finiteness of 
$\nu\otimes Q$ and $Q_s(\xi=0)=0$ for $\nu$-a.e. $s\in\bS,$ have been 
shown in Proposition \ref{CharPalm}. Equation \eqref{MeckeGl} is a special 
case of \eqref{Neveu2}.

Conversely assume the regularity conditions and 
\eqref{MeckeGl}. By Proposition \ref{CharPalm} it is enough to show 
that this implies \eqref{Sym}. By means of \eqref{OrbitRep} we have
\begin{align*}
\iiint f(\omega,s,t)\xi(\omega,dt)Q_s(d\omega)\nu(ds&) =
\iint\BE_{Q_s}\int f(\theta_e,s,t)\xi(dt)\mu_b(ds)\nu^*(db) \\
&= \iint\BE_{Q_{hb}}\int f(\theta_e,hb,t)\xi(dt)
\lambda(dh)\nu^*(db)\\
&= \iint\BE_{Q_{b}}\int f(\theta_h,hb,ht)\xi(dt)
\lambda(dh)\nu^*(db),
\end{align*}
where we used invariance of $Q$ and $\xi$ in the last step. 
Using the stochastic kernel $\kappa$ this last expression can be 
written as 
\begin{align*}
\int\BE_{Q_{b}}\iiint f(\theta_h,hb,hg\beta(t))\kappa_{\beta(t),t}(dg)
\xi(dt)\lambda(dh)\nu^*(db),
\end{align*}
and this equals
\begin{align*}
\int\BE_{Q_{b}}\iiint f(\theta_{hg^{-1}},hg^{-1}b,h\beta(t))
\lambda(dh)\Delta(g^{-1})\kappa_{\beta(t),t}(dg)\xi(dt)\nu^*(db),
\end{align*}
by Fubini's theorem and a characteristic property of the modular function. 
Now apply \eqref{MeckeGl} to the function 
$(\omega,s,t)\mapsto\int f(\theta_{h}\omega,hs,ht)\lambda(dh)$ 
to write this as
\begin{align*}
\int\BE_{Q_{b}}\iint f(\theta_{h},ht,hb)\lambda(dh)\xi(dt)\nu^*(db).
\end{align*}
By Fubini's theorem and invariance of $Q$ and $\xi$ this can be written as
\begin{align*}
\iint\BE_{Q_{hb}}\int f(\theta_{e},t,hb)\xi(dt)\lambda(dh)&\nu^*(db) \\
&= \iint\BE_{Q_{s}}\int f(\theta_{e},t,s)\xi(dt)\mu_b(ds)\nu^*(db)\\
&=\iiint f(\omega,t,s)\xi(\omega,dt)Q_s(d\omega)\nu(ds),
\end{align*}
where we have used \eqref{OrbitRep} for the second equality.\qed

\section{The mass-transport principle}\label{sectransport}

In this section we will show that Theorem \ref{thlastcycle} 
contains a mass-conservation law. 
Recall our basic properness assumption for the operation 
of $\bG$ on $\bS$. Let $M$ denote a $\sigma$-finite invariant 
measure on $\bS\times\bS$, 
which is given the interpretation that $M(C\times D)$ 
represents mass being transported out of $C$ into $D$. Then the main 
result of Subsection \ref{subdet}
says that on any set $B$ with a symmetry property with respect to 
the operating group $\bG$ the mass transported out 
of $B$ - suitably weighted in the non-unimodular case - equals the 
total mass transported 
into $B$. For a precise formulation of this
{\em mass-transport principle} (short: {\em MTP}) 
fix $k>0$ as in \eqref{KallFunk} and define
a measurable function $\tilde\Delta:\bS\times\bS\to(0,\infty)$ by 
\begin{align}\label{deltak}
\tilde\Delta(s,t)
:=\frac{\mu_t k}{\mu_{\beta(t)}k}\frac{\mu_{\beta(s)}k}{\mu_s k}.
\end{align}
Note that $\tilde\Delta(s,gs)=\Delta(g^{-1})$, $g\in\bG$, $s\in\bS$. 
A {\em symmetric set} is a set $B\in\cS$ satisfying 
$0<\mu_b(B)=\mu_{b'}(B)<\infty$, $b,b'\in O$.  The following
result will be obtained as a special case
of Theorem \ref{shortMTP} in Subsection \ref{subdet}.

\begin{theorem}\label{shortMTPonSets}
For any $\sigma$-finite and jointly invariant measure $M$ on $\bS\times\bS$ 
and any symmetric $B\in\mathcal{S}$ we have
\begin{align*}
\int \tilde\Delta(s,t)\I_B(s) M(d(s,t)) = \int \I_B(t)M(d(s,t)).
\end{align*}
\end{theorem}

In the transitive case, any set $B\in\cS$ with
positive and finite invariant measure is symmetric 
and Theorem \ref{shortMTPonSets} can be simplified as follows.

\begin{corollary}
If $\bG$ operates transitively on $\bS$ then for any $\sigma$-finite 
invariant measure $M$ on $\bS\times\bS$ and any $B\in\cS$ we have
\begin{align*}
\int \tilde\Delta(s,t)\I_B(s) M(d(s,t)) = \int \I_B(t)M(d(s,t)),
\end{align*}
where $\tilde\Delta(s,t) = \Delta(g^{-1})$ for any $g\in\bG$ with $gs=t$.
\end{corollary}

Up to an integrability issue Theorem \ref{shortMTPonSets}
implies the following stochastic analogue for 
invariant random measures. In  
Subsection \ref{subtrpri} we shall derive this result 
from the transport formula of Theorem \ref{thlastcycle}. 
A function on $\Omega\times \bS$ (or other
product spaces) is {\em invariant} if it is invariant under
joint shifts of the arguments.

\begin{theorem}\label{MTP}
Let $\xi,\eta$ be invariant random measures 
on $\bS$, and $\gamma,\delta$ invariant kernels 
from $\Omega\times\bS$ to $\bS$ such 
that \eqref{lastT*} holds for $\BP$-a.e.\ $\omega$.
Then for a symmetric set $B\in\mathcal{S}$ and any
invariant $h\in(\cA\otimes\cS\otimes\cS)_+$ we have
\begin{align}\label{set-form}
\BE\iint \tilde\Delta(s,t)\I_B(s)h(\theta_e,s,t)&\gamma(s,dt)\xi(ds) \notag\\
& = \BE\iint \I_B(t)h(\theta_e,s,t)\delta(t,ds)\eta(dt).
\end{align}
\end{theorem}

\subsection{Deterministic transport principle}\label{subdet}

Imagine there is some mass distributed over 
the space $\bS$ and that we transport from each $s\in\bS$ to each $t\in\bS$ some 
mass $m(s,t)$ in an invariant way, i.e.\ we assume that 
$$
m(gs,gt)=m(s,t),\quad g\in\bG.
$$ 
One  might guess that the total mass 
being transported from one orbit in $\bS$ to some fixed point in $\bS$ does 
only depend on the orbit of this fixed point and in fact equals the 
total mass being transported from any representative of the initial 
orbit into the whole orbit of our target point. 

To make this precise fix some representatives $b,b'\in O$ and consider 
their corresponding orbital invariant measures $\mu_b, \mu_{b'}$ 
which are $\sigma$-finite by our properness assumption. 
Recall the definition \eqref{lastDelta*} of $\Delta^*$.
The calculation
\begin{align*}
\int m(b,s)\Delta^*(s)\mu_{b'}(ds) &=\int m(b,gb')\Delta(g^{-1})\lambda(dg) = 
\int m(b,g^{-1}b')\lambda(dg) \\
&=\int m(gb,b')\lambda(dg) = \int m(s,b')\mu_{b}(ds)
\end{align*}
yields a basic balance equation between any two orbits:

\begin{corollary}\label{orbitbalance}
Let $m\in(\cS\otimes\cS)_+$ be invariant. Then 
\begin{align}
\label{elemdetmtp}
\int m(b,s)\Delta^*(s)\mu_{b'}(ds) = \int m(s,b')\mu_b(ds),
\quad b,b'\in O.
\end{align}
\end{corollary}

Corollary \ref{orbitbalance} is in fact a consequence of
\eqref{Neveu}. To see this 
specialize $\xi$ and $\eta$ to $\xi:=\mu_b$, $\eta:=\mu_{b'}$, choose respective 
supporting measures $\nu_{\mu_b}:=\mu_b$ and $\nu_{\mu_{b'}}:=\mu_{b'}$, 
note that
$Q_{\xi,b}=\BP=Q_{\eta,b}$, $b\in O$, and that $\mu^*_b=\delta_b$, replace
$f(\omega,g,s,t):=g(\omega)m(s,t)$ where $g$ is chosen such that 
$\BE_\BP[g]=1$, factor
out the $\sigma$-finite and invariant $\BP$ on both sides and finally use
invariance of $m$.

Corollary \ref{orbitbalance} leads to the following 
mass-transport principle on any system $O$ of orbit representatives. 
Given two invariant measures $\mu$ and $\nu$ on $\bS$ we may interpret
$\mu$ as  mass distributed within $\bS$ while $\nu$ 
on the other hand can be thought of holes where mass can be stored. 
Consider invariant  
kernels $\gamma$ and $\delta$ on $\bS$, where an application of 
$\gamma$ to $\mu$ may 
be thought of resizing and redistributing the mass $\mu$, while an application of 
$\delta$ to $\nu$ may be interpreted as streching or shrinking and relocating 
the holes $\nu$. Therefore one might call $\gamma$ and $\delta$
{\em weighted transport kernels}, see \cite{LaTh08}. 

\begin{corollary}\label{Detmtprep}
Let $\mu,\nu$ be $\sigma$-finite invariant measures on $\bS$ and 
$\gamma,\delta$ invariant kernels on $\bS$ satisfying
\begin{align}\label{**}
\iint\I\{(s,t)\in\cdot\}\gamma(s,dt)\mu(ds)
=\iint\I\{(s,t)\in\cdot\}\delta(t,ds)\nu(dt).
\end{align} 
Then for any invariant $m\in(\cS\otimes\cS)_+$ 
\begin{align}\label{detmtprep}
\iint m(b,t)\Delta^*(t)\gamma(b,dt)\mu^*(db) = \iint m(s,b)\delta(b,ds)\nu^*(db).
\end{align}
\end{corollary}
{\em Proof.}
Specializing $\xi:=\mu$, $\eta:=\nu$ in Theorem $\ref{thlastcycle}$
and following similar steps as above yields the result.  \qed

\vspace{0.3cm}

The following version of the above deterministic MTP has been 
established in \cite{B:L:P:S:99} in the case 
of finitely many orbits (see also \cite{LP:book}, Chapter 8).

\begin{corollary} 
Let $\bG$ operate properly on the countable set $\bS$.
Then we have for invariant $m\in(\cS\otimes\cS)_+$
\begin{align*}
\sum_{b,b'\in O}\sum_{s\in\bG b'} m(b,s) \lambda(\bG_{s,s})= 
\sum_{b,b'\in O}\sum_{s\in\bG b'} \lambda(\bG_{b',b'})m(s,b)
\end{align*}
and if $\bG$ is unimodular even
\begin{align*}
\sum_{b\in O}\frac{1}{\lambda(\bG_{b,b})}\sum_{s\in\bS} m(b,s) = 
\sum_{b\in O}\frac{1}{\lambda(\bG_{b,b})}\sum_{s\in\bS} m(s,b).
\end{align*}
\end{corollary}
{\em Proof.}
Putting $\mu=\nu:=\sum_{b\in O}\mu_b$ (hence $\mu^*=\nu^*=\sum_{b\in O}\delta_b$) 
and $\gamma(s,\cdot) \equiv \delta(t,\cdot) \equiv \sum_{b\in O}\mu_b$ in \eqref{detmtprep} 
yields after a similar step as in \eqref{intmu}
\begin{align*}
\sum_{b,b'\in O}\lambda(\bG_{b',b'})\sum_{s\in\bG b'} m(b,s)\Delta^*(s) = 
\sum_{b,b'\in O}\lambda(\bG_{b',b'})\sum_{s\in\bG b'} m(s,b)
\end{align*}
for any measurable invariant $m$. Using \eqref{quothaar} this simplifies to
\begin{align*}
\sum_{b,b'\in O}\sum_{s\in\bG b'} m(b,s) \lambda(\bG_{s,s})= 
\sum_{b,b'\in O}\sum_{s\in\bG b'} \lambda(\bG_{b',b'})m(s,b).
\end{align*}
If $\bG$ is unimodular then $\lambda(\bG_{b',b'})$ on the right may be
replaced by $\lambda(\bG_{s,s})$ and we may also
replace $m$ by the $\bG$-invariant function $m(t,s)\frac{1}{\lambda(\bG_{t,t})
\lambda(\bG_{s,s})}$ which then yields the assertion.\qed

\vspace{0.3cm}

We seek for an appropriate formulation of 
Corollary \ref{Detmtprep} without use of a fixed 
system of representatives of the orbits. For this choose some $v\in\cS_+$ 
with the property that $0<\mu_sv<\infty$ for each $s\in\bS$ and define
similarly as in \eqref{deltak}
\begin{align}\label{Delta}
\Delta_v(s,t):=\frac{\mu_tv}{\mu_sv},\quad s,t\in\bG.
\end{align} 
By means of another function $w\in\cS_+$ with $0<\mu_sw<\infty$ for any $s\in\bS$
being compatible with $v$ in the sense that
\begin{align}\label{balance}
\frac{\mu_bw}{\mu_bv}=\frac{\mu_{b'}w}{\mu_{b'}v},\quad b,b'\in O,
\end{align}
we may, since $\Delta^*(s)=\mu_s(v)/\mu_{\beta(s)}(v)$, express $\Delta^*$ as
\begin{align}\label{identity}
\Delta^*(s) = \Delta^*(s)\frac{\mu_{\beta(s)}(v)}{\mu_{\beta(s)}(w)}
\frac{\mu_b(w)}{\mu_{b}(v)} = \Delta_v(b,s)\frac{\mu_b(w)}{\mu_{\beta(s)}(w)},
\quad b\in O, s\in\bS,
\end{align}
where we have used assumption \eqref{balance}. 
 
\begin{proposition}\label{mtpinttransportkernels}
Let $\mu$, $\nu$ be $\sigma$-finite invariant measures on $\bS$ and 
$\gamma$ and $\delta$ be invariant kernels on $\bS$ satisfying \eqref{**}.
Let $v,w\in\cS_+$ be as in \eqref{balance} and $m\in(\cS\otimes\cS)_+$ 
be invariant. Then 
\begin{align*}
\iint \Delta_v(s,t)w(s)m(s,t)\gamma(s,dt)\mu(ds) = 
\iint w(t)m(s,t)\delta(t,ds)\nu(dt).
\end{align*}
\end{proposition}
{\it Proof.} Replacing $\Delta^*(t)$ in \eqref{detmtprep} by the left-hand side 
of \eqref{identity}, we get
\begin{align*}
\iint m(b,t)\Delta_v(b,t)\frac{\mu_b(w)}{\mu_{\beta(t)}(w)}\gamma(b,dt)\mu^*(db)
=\iint m(s,b)\delta(b,ds)\nu^*(db).
\end{align*}
Applying this with $m$ replaced by the
invariant function $m(s,t)\mu_{\beta(t)}(w)$ yields
\begin{align*}
\iint m(b,t)\Delta_v(b,t)\mu_b(w)\gamma(b,dt)\mu^*(db) = 
\iint m(s,b)\mu_b(w)\delta(b,ds)\nu^*(db).
\end{align*}
Fubini's theorem yields
\begin{align*}
\iiint m(b,t)\Delta_v(b,t)w(s)&\gamma(b,dt)\mu_b(ds)\mu^*(db) \notag\\
&= \iiint m(s,b)w(t)\delta(b,ds)\mu_b(dt)\nu^*(db).
\end{align*}
By invariance of $m\Delta_v$ and invariance of $\gamma$ 
the integral
$\int m(s,t)\Delta_v(s,t)\gamma(s,dt)$ does not depend on $s\in Gb$ for a fixed $b\in O$.
A similar remark applies to the inner integral
of the above right-hand side. Applying \eqref{OrbitRep} we obtain
the asserted identity. \qed

\vspace{0.3cm}

Using a result of Kallenberg in \cite{Kall07}, 
Proposition \ref{mtpinttransportkernels} 
can be formulated equivalently in a shorter way.

\begin{theorem}\label{shortMTP}
Let $v,w\in\cS_+$ satisfy \eqref{balance}. Then for any jointly 
invariant $\sigma$-finite measure $M$ on $\bS\times\bS$ we have
\begin{align}
\int \Delta_v(s,t)w(s)M(d(s,t)) = \int w(t)M(d(s,t)).
\end{align}
\end{theorem}
{\em Proof.}
According to Kallenberg \cite{Kall07} Corollary 3.6, there exist 
invariant disintegrations of the form
\begin{align*}
M(d(s,t))=\gamma(s,dt)\mu(ds), \quad M(d(s,t))=\delta(t,ds)\nu(dt),
\end{align*}
where $\mu,\nu$ are invariant $\sigma$-finite measures on $S$ and
$\gamma,\delta$ are invariant kernels on $S$.
Then \eqref{**} holds, and
Proposition \ref{mtpinttransportkernels} implies the assertion.\qed

\vspace{0.3cm}

{\em\noindent Proof of Theorem \ref{shortMTPonSets}:}
Just note that 
$v(s):=k(s)/\mu_{\beta(s)}k$ and $w:=\I_B$ for a symmetric $B\in\mathcal{S}$ 
is an admissable choice in Theorem \ref{shortMTP} since 
\begin{align*}\tag*{\qed}
\frac{\mu_b v}{\mu_b(B)} = \frac{1}{\mu_b(B)} = \frac{1}{\mu_{b'}(B)} 
= \frac{\mu_{b'} v}{\mu_{b'}(B)},\quad b,b'\in O. 
\end{align*}

\subsection{Transport principle for stationary random measures}\label{subtrpri}

In this section we show how Theorem \ref{MTP} follows from the transport formula
in Theorem \ref{thlastcycle}. To this end we begin with a 
modification of the refined Campbell theorem for 
invariant functions.

\begin{lemma}
Let $\xi$ be an invariant random measure on $\bS$ and $(\nu_\xi,Q_\xi)$
an invariant Palm pair of $\xi$. Then for any $\bG$-invariant 
$f\in(\cA\otimes\cS)_+$ and any $v\in\cS_+$
\begin{align}\label{inversion}
\int \BE_{Q_{\xi,b}}[f(\theta_e,b)]\mu_b(v)\nu_\xi^*(db) = 
\BE_\BP\int f(\theta_e,s)v(s)\xi(ds).
\end{align}
\end{lemma}
{\em Proof.}
The left side equals 
$$
\BE\iint f(\theta_g^{-1},\beta(s))v(g\beta(s))\kappa_{\beta(s),s}(dg)\xi(ds)
$$
by Proposition \ref{refC}. 
Invariance of $f$ yields the result. \qed

\vspace{0.3cm}

Applying this modified Campbell formula for invariant 
functions to the transport formula
of Theorem \ref{thlastcycle} yields Theorem \ref{MTP}:

\vspace{0.3cm}

{\em Proof of Theorem \ref{MTP}.}
Replacing $f(\omega,g,s,t)$ 
in \eqref{Neveutrans} by $h(\omega, gs,t)$ for 
$h\in(\cA\otimes\cS\otimes\cS)_+$ yields 
\begin{align*}
\int\BE_{Q_{\eta,b}}\iint h(\theta_g^{-1}, &g^{-1}b,\beta(s))
\Delta^*(s)\kappa_{\beta(s),s}(dg)\delta(b,ds)\nu^*_\eta(db) \notag \\
 &=\int\BE_{Q_{\xi,b}}\iint h(\theta_e, g\beta(s),b)
\kappa_{\beta(s),s}(dg)\gamma(b,ds)\nu^*_\xi(db).
\end{align*}
If $h$ is invariant this reduces to
\begin{align*}
\int\BE_{Q_{\eta,b}}\int h(\theta_e, b,s)\Delta^*(s)&\delta(b,ds)\nu^*_\eta(db) = 
\int\BE_{Q_{\xi,b}}\int h(\theta_e, s,b)\gamma(b,ds)\nu^*_\xi(db).
\end{align*}
Substituting $\Delta^*(s)$ via (\ref{identity}) and replacing
$h(\omega,t,s)$ by $h(\omega,t,s)\mu_{\beta(s)}(w)$ 
yields for any invariant $h$
\begin{align*}
\int\BE_{Q_{\eta,b}}\int h(\theta_e,b,s)\Delta_v(b,s)\mu_b(w)
&\delta(b,ds)\nu^*_\eta(db) \notag \\
&= \int\BE_{Q_{\xi,b}}\int h(\theta_e,t,b)\mu_b(w)\gamma(b,dt)\nu^*_\xi(db).
\end{align*}
Here $(\omega,t)\mapsto\int h(\omega,t,s)\Delta_v(t,s) \delta(\omega,t,ds)$ 
and 
$(\omega,s)\mapsto\int h(\omega,t,s)\gamma(\omega,s,dt)$ are clearly 
invariant by invariance of the transports $\gamma,\delta$ and 
the function $h$. Thus 
\eqref{inversion} may be applied to both sides of the last equation to yield 
\begin{align*}
\BE\iint h(\theta_e,t,s)\Delta_v(t,s)w(t)\delta(t,ds)\eta(dt)=
\BE\iint h(\theta_e,t,s)w(s)\gamma(s,dt)\xi(ds).
\end{align*}
Since $\Delta_v(s,t) = \Delta_v(t,s)^{-1}$ this is clearly equivalent to 
saying that for any invariant $h$
\begin{align}\label{6.12}
\BE\iint h(\theta_e,t,s)w(t)\delta(t,ds)&\eta(dt) \notag \\
&=\BE\iint h(\theta_e,t,s)\Delta_v(s,t)w(s)\gamma(s,dt)\xi(ds).
\end{align}
Choosing $v,w$ as in the proof of Theorem \ref{shortMTPonSets}  yields the result. \qed

\vspace{0.3cm}

Putting in \eqref{6.12} $h(\omega,t,s)=m(\omega,s,t)$ for some invariant  
$m\in(\mathcal{A}\otimes\cS\otimes\cS)_+$,
$\gamma(\omega,s,dt):=\eta(\omega,dt)$ and
$\delta(\omega,t,ds):=\xi(\omega,ds)$  we obtain
the following corollary.

\begin{corollary}
For $v,w\in\cS_+$ satisfying \eqref{balance}, invariant random 
measures $\xi,\eta$ on $\bS$ and invariant 
$m\in(\mathcal{A}\otimes\cS\otimes\cS)_+$ we have
\begin{align}\label{mtponrm}
\BE\iint \Delta_v(s,t)w(s)m(\theta_e,s,t)&\xi(ds)\eta(dt) \notag \\
&= \BE\iint w(t)m(\theta_e,s,t)\xi(ds)\eta(dt).
\end{align}
\end{corollary}

\vspace{0.3cm}

\noindent
{\bf Note:} Some of the main results of this paper were
presented in the workshop ``New Perspectives in Stochastic Geometry''
at the Mathematisches Forsch\-ungs\-institut Oberwolfach,
October 2008.

\noindent
{\bf Note:} While this paper was under review, we became
aware of the preprint \cite{Ka09} by Olav Kallenberg,
written independently of this one and finished
at about the same time.
Among many other things the author proves there
the existence of a universally measurable choice function under the sole assumption
of properness as well as Theorem \ref{charMecke}.

\end{document}